\documentclass[11pt]{article}

\usepackage{amsmath,dsfont,mathrsfs}
\usepackage{verbatim}
\usepackage{amssymb}
\usepackage{latexsym}
\usepackage[latin1]{inputenc}
\usepackage[T1]{fontenc}
\usepackage{pict2e}
\usepackage{graphicx}

\newtheorem{Theorem}{Theorem}[section]

\newtheorem{Lemma}{Lemma}[section]

\newenvironment{Proof}{\vspace{1.5ex}{\sc Proof}.}{\vspace{2ex}}

\def\1e{\mathds{1}}

\def\Ne{\mathds{N}}

\def\Re{\mathds{R}}

\def\vgc{\mathbf{c}}

\def\vgf{\mathbf{f}}
\def\vgg{\mathbf{g}}

\def\vgx{\mathbf{x}}
\def\vgy{\mathbf{y}}

\def\vg0{\mathbf{0}}

\def\E{\mathscr{E}}

\def\P{\mathscr{P}}

\def\R{\mathscr{R}}

\def\thgg{\boldsymbol{\theta}}

\def\xgg{\boldsymbol{\xi}}

\def\d{\displaystyle}

\def\eop{\rule{0.7ex}{0.7ex}}
\def\eqd{:=}

\def\segoc#1#2{(#1, #2]}

\def\scal#1#2{
    \left\langle #1,#2\right\rangle
    }
    
\def\scalf#1#2{
    \left\langle #1,#2\right\rangle_F
    }
    
\def\set#1#2{
    \left\{#1|\vphantom{#1#2}#2\right\}
    }

\def\norm#1{\left\|#1\right\|}

\def\normf#1{\left\|#1\right\|_F}

\def\argmin#1{\mathop{\mathrm{argmin}}#1}

\def\tr#1{\mathop{\mathrm{tr}}#1}

\def\ran#1{\mathop{\mathrm{ran}}#1}
\def\ker#1{\mathop{\mathrm{ker}}#1}

\def\dom#1{\mathop{\mathrm{dom}}#1}

\def\sgntile#1{\Gamma}
\def\htile#1{\Gamma_e}

\def\ud{\,\mathrm{d}}

\begin{document}


\title{Tomography by Fourier synthesis}
\author{P. Mar\'echal and A. Elasmai\\[.5ex]
Institut de Math\'ematiques de Toulouse}
\date{\today}
\maketitle



%

%


%
\abstract{
We consider a particular approach to the regularization of the inverse
problem of computerized tomography. This approach is based on notions
pertaining to Fourier synthesis. It refines previous contributions,
in which the preprocessing of the data was
performed according to the Fourier slice theorem. Since real models
must account for the geometrical system response and possibly
Compton scattering and attenuation, the Fourier slice theorem
does not apply, yielding redefinition of the preprocessing.
In general, the latter is not explicit, and must be performed
numerically. The most natural choice of preprocessing involves
the computation of unstable solutions. A proximal strategy is
proposed for this step, which allows for accurate
computations and preserves global stability of the reconstruction process.
}

                                        


\vspace{2ex}

{\bf Warning.}
This paper was submitted to the Journal of Inverse and Ill-Posed Problems on May 10, 2010,
accepted there on November 28, 2010, and however never published.

\section{Introduction}

From the mathematical viewpoint, tomography consists in solving
an ill-posed operator equation of the form $Rf_0=g$, where $f_0$
is the unknown image, $R$ is a linear operator (a simplified version
of which being the Radon transformation) and $g$ is the data
provided by the imaging device ({\it e.g.} a SPECT or PET camera).
The data~$g$ is usually referred to as the {\sl sinogram}.

In the ideal model, $R$ is the standard Radon operator, given by
\begin{equation}
\label{RadonOperator}
(Rf)(\thgg,r)=\int f(\vgx)\delta(r-\scal{\thgg}{\vgx})\ud\vgx.
\end{equation}
Here, $\thgg=(\cos\phi,\sin\phi)$ is an element of the unit circle
in $\Re^2$, $r$~is a real variable and $\delta$ denotes the
Dirac {\sl delta function}, so that the above integral is in fact the
(one-dimensional) integral of~$f$ over the line determined by
the equation $r=\scal{\thgg}{\vgx}$.
In this case, the celebrated {\sl Fourier Slice Theorem}
allows for considering the constraints on~$f_0$ imposed by the
equation $Rf_0=g$ as constraints on the Fourier transform of~$f_0$,
these constraints appearing on a radial domain of the Fourier plane.
This observation led the authors of~\cite{m-togane-c,dmg-m-et-al}
to regularize the ill-posed operator equation by means of concepts
pertaining to {\sl Fourier synthesis}~\cite{lannes-rc}.

Recall that {\sl Fourier synthesis} refers to the generic problem of recovering
a function~$f_0$ from a partial and approximate knowledge of its Fourier transform.
In \cite{lannes-rc}, the analysis of spectral properties of the truncated
Fourier operator led the authors to introduce a regularization principle
which can be regarded as a reformulation into a well-posed problem
of Fourier interpolation. The original problem of recovering the unknown object
is replaced by that of recovering a limited resolution
version of it, namely, $\phi_\beta\ast f_0$, where $\phi_\beta$
is some convolution kernel (or point spread function).
In~\cite{almasa}, inspired by well-known results from the approximation of $L^p$-functions
by means of mollifiers, the authors established variational results
for this type of regularization. They regarded~$\phi_\beta$ as a member of the
one-parameter family $\set{\phi_\beta}{\beta\in\Re_+^*}$ defined by
$$
\phi_\beta(\vgx)\eqd\frac{1}{\beta^2}\phi\left(\frac{\vgx}{\beta}\right),
$$
and studied the behavior of the reconstructed object as $\beta$ goes to zero.

Notice that mollifiers were also considered in~\cite{louis-maass}. By writing
the {\sl mollified} function $f_\beta\eqd\phi_\beta\ast f$ as
$$
f_\beta(\vgx)\eqd\int\phi_\beta(\vgx-\vgy)f(\vgy)\ud\vgy=\scal{\varphi_\beta(\vgx,\cdot)}{f},
$$
with $\varphi_\beta(\vgx,\vgy)\eqd\phi_\beta(\vgx-\vgy)$, the
duality associated to the underlying inner product is used to approximate $f_\beta$ {\it via}
$$
f_\beta(\vgx)\simeq\scal{R^*\psi_\beta(\vgx)}{f}=\scal{\psi_\beta(\vgx)}{Rf},
$$
in which $\psi_\beta(\vgx)\simeq (R^*)^{-1}\varphi_\beta(\vgx,\cdot)$. This approach
is referred to as the method of {\sl approximate inverses}.

Our approach inherits features from both Fourier synthesis and the approximate inverses.
In Section~\ref{methodology}, we give a precise description of our methodology in a discrete
setting. We shall focus on the problem of Computed Tomography.
We shall introduce the notion of
{\sl pseudo-commutant} of a matrix with respect to another.
In section~\ref{computational},
we discuss various computational strategies, especially concerning what we call
{\sl preprocessing of the data}, which consists in applying the aformentioned
pseudo-commutant. Finally, in Section~\ref{numerical}, we prove the numerical feasibility
of our approach by means of reconstructions in emission tomography, in which the
evolution towards realistic {\sl projectors} made it necessary to give up
the mere application of the Fourier Slice Theorem.
Our numerical experiments will clearly indicate improvement
in terms of stability and image quality.

\section{A reconstruction methodology}
\label{methodology}

As outlined in the introduction, our aim is to reconstruct a smoothed version
$\phi_\beta\ast f_0$ of the original object~$f_0$. Recall that, if the operator
modelling data acquisition is actually the Radon operator~$R$ defined in
Equation~\eqref{RadonOperator}, it is easy to generate the data corresponding
to $\phi_\beta\ast f_0$. As a matter of fact, it results immediately from
the Fourier Slice Theorem that
\begin{equation}
\label{radon-convolution}
R(\phi_\beta\ast f_0)= R\phi_\beta\circledast Rf_0,
\end{equation}
where $\circledast$ denotes the convolution with respect to~$r$.
If~$g$ is an approximation of $Rf_0$, then $R\phi_\beta\circledast g$ will be
an approximation of $R(\phi_\beta\ast f_0)$. In terms of operators, this amounts to
the existence of an operator $\Phi_\beta$ such that $RC_\beta=\Phi_\beta R$,
where $C_\beta$ is the convolution operator:
$C_\beta f\eqd\phi_\beta\ast f$.
Now, realistic models describing the data acquisition in emission or transmission
tomography involve operators
which are not the exact Radon transformation, and which do not satisfy,
in general, Equation~\eqref{radon-convolution}.

In a recent paper~\cite{xapi}, an extension of the regularization by mollification
was proposed, in order to cope with operators~$R$ for which it is not possible to
explicitly find an operator $\Phi_\beta$ such that $RC_\beta=\Phi_\beta R$.
The idea consists in defining $\Phi_\beta$ as an operator minimizing
$$
X\mapsto\norm{RC_\beta-XR},
$$
where $\norm{\cdot}$ is some operator norm. In~\cite{xapi}, the authors
focus on the infinite dimensional setting, while considering general operators.
Here, we concentrate on the discrete case. Therefore, from now on,
$R$, $C_\beta$ and $X$ are matrices of respective sizes $m\times n$,
$n\times n$ and $m\times m$, and we shall deal with the minimization
of the Frobenius norm the matrix $RC_\beta-XR$ (with respect to~$X$).

Recall that the Frobenius norm of a matrix $M\in\Re^{m\times n}$,
denoted by $\normf{M}$, is the Euclidean norm of~$M$
regarded as a vector in $\Re^{mn}$, and that the corresponding inner product,
denoted by $\scalf{\cdot}{\cdot}$, satisfies:
$$
\forall M,N\in\Re^{m\times n},\quad
\scalf{M}{N}=\tr{(MN^\top)}=\tr{(N^\top M)}.
$$
The following theorem can be found in various sources and various forms.
We provide a proof for the sake of completeness.
Recall that a matrix~$X\in\Re^{n\times m}$ is the pseudo-inverse $M^\dagger$
of $M$ if and only if it satisfies $MXM=M$, $XMX=X$, $(MX)^\top=MX$
and $(XM)^\top=XM$.

\begin{Theorem}\sf
Let $A$ and $R$ be real matrices of size $m\times n$.
Then the matrix $AR^\dagger$ minimizes the function
$f(X)\eqd\normf{A-XR}^2/2$ over $\Re^{m\times m}$.
Among all minimizers, it is the one with minimum Frobenius norm.
\end{Theorem}

\Proof{
Clearly, $f$ is convex, indefinitely differentiable, and
its gradient at~$X$ is equal to $(XR-A)R^\top$. Therefore,
$$
\nabla f(AR^\dagger)=
A(R^\dagger RR^\top-R^\top)=
A\big((RR^\dagger R)^\top-R^\top\big)=0.
$$
This proves that $AR^\dagger$ minimizes~$f$.
Furthermore, since the Frobenius norm is strictly convex, minimizers of~$f$
can only differ by matrices in the kernel of the linear mapping
$\R\colon X\mapsto XR$ (from $\Re^{m\times m}$ to~$\Re^{m\times n}$).
Now, for all $K\in\Re^{m\times m}$,
$$
\scalf{AR^\dagger}{K}=
\tr{\big(AR^\dagger RR^\dagger K^\top\big)}=
\tr{\big(AR^\dagger(R^\dagger)^\top(KR)^\top\big)}.
$$
This implies that every matrix $K$ in the kernel of~$\R$
is orthogonal (for $\scalf{\cdot}{\cdot}$) to $AR^\dagger$.
The desired conclusion follows by Pythagoras' theorem.~\eop}

Now, letting $A$ be the matrix $RC_\beta$, we see that
the minimum Frobenius norm minimizer of the matrix
function $X\mapsto\normf{RC_\beta-XR}$ is the matrix
$$
\Phi_\beta\eqd RC_\beta R^\dagger\in\Re^{m\times m}.
$$
In order to be as consistent as possible with our aim, that is,
with the reconstruction of the convolution of the original image by
our point spread function~$\phi_\beta$, we must replace
the original sinogram with its transformation by $RC_\beta R^\dagger$.

Let $\vgf$ and $\vgg$ denote the discrete versions of~$f$ and~$g$,
respectively. Consider the decomposition of the generic image~$\vgf$ as the sum
of its low and high frequency components:
$$
\vgf=C_\beta \vgf+(I-C_\beta)\vgf.
$$
We define the reconstructed image as the solution to the following
optimization problem:
$$
(\P)\quad\left|
\begin{array}{rl}
\hbox{Minimize}&\d\frac{1}{2}\norm{RC_\beta R^\dagger \vgg-R\vgf}^2+
\frac{\alpha}{2}\norm{H_\beta \vgf}^2\\[1.5ex]
\hbox{subject to}& \vgf\geq\vg0,
\end{array}\right.
$$
in which $\alpha$ is a positive weight and $H_\beta\eqd I-C_\beta$.
From the computational viewpoint, the difficult part is the estimation
of the {\sl regularized data} $RC_\beta R^\dagger \vgg$.
The reason is, of course, the ill-posedness of the model matrix~$R$,
which yields ill-posedness of the computation of $R^\dagger \vgg$.
This difficulty will be addressed in the next section.

\section{Computational aspects}
\label{computational} 

In this section, we address the computation of the regularized
data, that is, of $RC_\beta R^\dagger\vgg$. The resolution of problems such as~$(\P)$
is quite standard, and will not be discussed here.

Due to the dimension of the involved
matrices, it seems unreasonable to actually compute either $RC_\beta R^\dagger$
or $R^\dagger$. The most natural strategy consists in computing $R^\dagger\vgg$
as the (unstable) minimum norm least square solution of the original system,
and then in applying $RC_\beta$ to the obtained solution. 

It may seem inappropriate to initialize a regularization scheme by computing
an unstable solution. It is important, however, to realize that ill-posedness encompasses
two aspects (which are related): the numerical inaccuracy induced by
the poor conditioning of the system, and the propagation of measurement errors~$\delta\vgg$.
We conjecture, at this point, that the former can be dealt with
by means of a {\sl proximal} strategy, while the latter is not crucial since the unstable
solution $R^\dagger\vgg$ is post-processed by the smoothing operator $C_\beta$
(and ultimately by~$R$). In essence, the proximal
point algorithm is a fixed point method. The reason for choosing
an iterative scheme which, incidentally, has the reputation of being slow,
is that it can give very accurate solutions. In fact, it may be used
for refining solutions provided by some other method.

In the next paragraph, we review a few aspects of the Proximal Point Algorithm (PPA),
applied to the computation of minimum norm least square solutions.

The Proximal Point Algorithm was introduced by Martinet~\cite{martinet} in 1970, in the context of variational
inequalities. It was then generalized by Rockafellar~\cite{rockafellar} to the
computation of zeros of maximal monotone operators (a particular case of which being the minimization
of a convex function).
In our context, that of computing the minimum norm least-square solution of the
linear system $R\vgf=\vgg$, the PPA consists in the following steps:
\begin{itemize}
\item[1.]
{\sc Initialization}: put $k=0$ and choose an initial point $\vgf_0\in\Re^n$.
\item[2.]
{\sc Iteration}:
form the sequence $(\vgf_k)_{k\geq 0}$ according to
$$
\vgf_{k+1}\eqd\argmin_{\vgf\in\Re^n}\left\{
\frac{1}{2}\norm{\vgg-R\vgf}^2+\frac{1}{2\lambda_k}\norm{\vgf-\vgf_k}^2\right\},
$$
where $(\lambda_k)_{k\geq 0}$ is a sequence of positive real numbers.
\end{itemize}
Observe that the function to be minimized in the iteration is
strictly convex, smooth and coercive, and differs from the Tikhonov
functional only by the subtraction of~$\vgf_k$ in the regularizing part.
It has the numerical stability inherent to Tikhonov regularization,
and may be computed by using any routine from either quadratic
optimization ({\it e.g.} the conjugate gradients) or linear algebra
(on the corresponding {\it regularized} normal equation).

In~\cite{rockafellar}, Rockafellar obtained convergence results
for the proximal point algorithm, in the general setting of
the computation of zeros of maximal monotone operators in real Hilbert
spaces. For the sake of simplicity, we restrict attention to the
minimization of any lower semi-continuous convex function on~$\Re^n$.
The proof of the following theorem is differed to the appendix.

\begin{Theorem}\sf
\label{theo-conv-prox}
Let $F\colon\Re^n\to\segoc{-\infty}{\infty}$ be a lower semi-continuous convex function
which is not identically equal to infinity and bounded below, and let $\vgx_0$ be any point in~$\Re^n$.
Consider the sequence $(\vgx_k)_{k\in\Ne}$ defined by
$$
\vgx_{k+1}=\argmin\set{F(\vgx)+\frac{1}{2\lambda_k}\norm{\vgx-\vgx_k}^2}{\vgx\in\Re^n},
$$
where $(\lambda_k)_{k\in\Ne}$ is a sequence of positive numbers. If the series
$\sum\lambda_k$ is divergent, then
$$
F(\vgx_k)\to\eta\eqd\inf\set{F(\vgx)}{\vgx\in\Re^n}
\quad\hbox{as}\quad
k\to\infty.
$$
If in addition the set of minimizers $S\eqd\set{\vgx\in\Re^n}{F(\vgx)=\eta}$ is non\-empty,
then the sequence $(\vgx_k)_{k\in\Ne}$ converges to a point~$\overline{\vgx}$ in~$S$.
\end{Theorem}

Now, let $F$ be the function $\vgf\mapsto\norm{\vgg-R\vgf}^2/2$. In this case, the set~$S$
of minimizers is the affine manifold
$$
S=\set{\vgf\in\Re^n}{R^\top R\vgf=R^\top\vgg}=\{R^\dagger\vgg\}+\ker{R}.
$$
Furthermore, writing the usual first order necessary optimality condition for~$\vgf_{k+1}$
yields
$$
\vgf_{k+1}-\vgf_k=-\lambda_k R^\top(R\vgf_{k+1}-\vgg)\in\ran{R^\top}=(\ker{R})^\perp.
$$
This shows that, if $\vgf_0=\vg0$, the limit point of the sequence $(\vgf_k)_{k\in\Ne}$
is nothing but $R^\dagger\vgg$ itself.

Notice that, in the case where $R$ is injective, the set of minimizers of
$\norm{R\vgf-\vgg}$ is reduced to the singleton $\{R^\dagger\vgg\}=\{(R^\top R)^{-1}R^\top\vgg\}$.
In this case, numerical errors in the proximal iteration are unimportant, since each new
iterate can be regarded as a new initial point. However, if~$R$ is not
injective, a component in the kernel of~$R$ may grow due to finite precision
in the iteration. In order to cope with this difficulty, the following strategy
may be adopted.
In practice, the smoothing properties of $RC_\beta$ are expected to
give rise to reasonable perturbation transmission: a small perturbation $\delta\vgf\in\Re^n$
will yield a reasonably small perturbation $RC_\beta\delta\vgf\in\Re^m$. Therefore, replacing~$R^\dagger$
with the usual Tikhonov approximation $(R^\top R+\varepsilon I)^{-1}R^\top$,
where~$I$ denotes the identity matrix of appropriate dimension and $\varepsilon$ is a
very small positive parameter, will give a
reasonable approximation of the regularized data:
$$
RC_\beta(R^\top R+\varepsilon I)^{-1}R^\top\vgg\approx
RC_\beta R^\dagger\vgg.
$$
In the above proximal scheme,
the iterations should then be replaced by 
$$
\vgf_{k+1}\eqd\argmin_{\vgf\in\Re^n}\left\{
\frac{1}{2}\norm{\vgg-R\vgf}^2+\frac{\varepsilon}{2}\norm{\vgf}^2+\frac{1}{2\lambda_k}\norm{\vgf-\vgf_k}^2\right\},
$$
which will clearly yield convergence to $(R^\top R+\varepsilon I)^{-1}R^\top\vgg$.

We stress that mathematical convergence is guaranteed whenever the series $\sum\lambda_k$
is divergent, and that this happens in particular if $\lambda_k\equiv\lambda$ with $\lambda>0$.
Notice that the smaller~$\lambda_k$, the smaller the size of the step $\vgf_{k+1}-\vgf_k$.
Consequently, it may be appropriate to start the algorithm with large values of~$\lambda_k$
yielding large (but poorly conditioned) steps, and then let $\lambda_k$ decrease
in order to have smaller but well-conditioned steps. 
Moreover, as indicated in~\cite{prox-inversion},
numerical accuracy in early iterations may be irrelevant: what really matters is the limit of the
proximal sequence.
We terminate this section with a few comments.

{\bf (1)}
If the size of~$R$ is not prohibitively large, the computation of~$R^\dagger$ may be performed
by means of the {\sl singular value decomposition} (SVD). The problems of tomographic reconstruction encountered
in practice make it difficult or even impossible to perform the SVD.
Approximate pseudo-inverses may then be obtained by truncating the SVD.
Such approximate solution may be used as initial points for the PPA.

{\bf (2)}
As observed in~\cite{prox-inversion}, the proximal iteration with constant sequence $\lambda_k\equiv\lambda$
belongs to the class of {\sl fixed point methods},
along with the algorithms of Jacobi, Gauss-Seidel, SOR and
SSOR. It is easy to check that $R^\dagger\vgg$ satisfies the fixed point equation
$\vgf=h(\vgf)$, with
$$
h(\vgf)\eqd B\vgf+\vgc,\quad
B\eqd(I+\lambda R^\top R)^{-1}
\quad\hbox{and}\quad
\vgc\eqd\lambda(I+\lambda R^\top R)^{-1}R^\top\vgg.
$$
Clearly, $h$ is a contraction and, if $R$ is injective, then $R^\top R$ is positive definite and $h$ is a strict
contraction.

\section{Numerical results}
\label{numerical}

\begin{figure}
\begin{center}
\includegraphics[width=3.7cm]{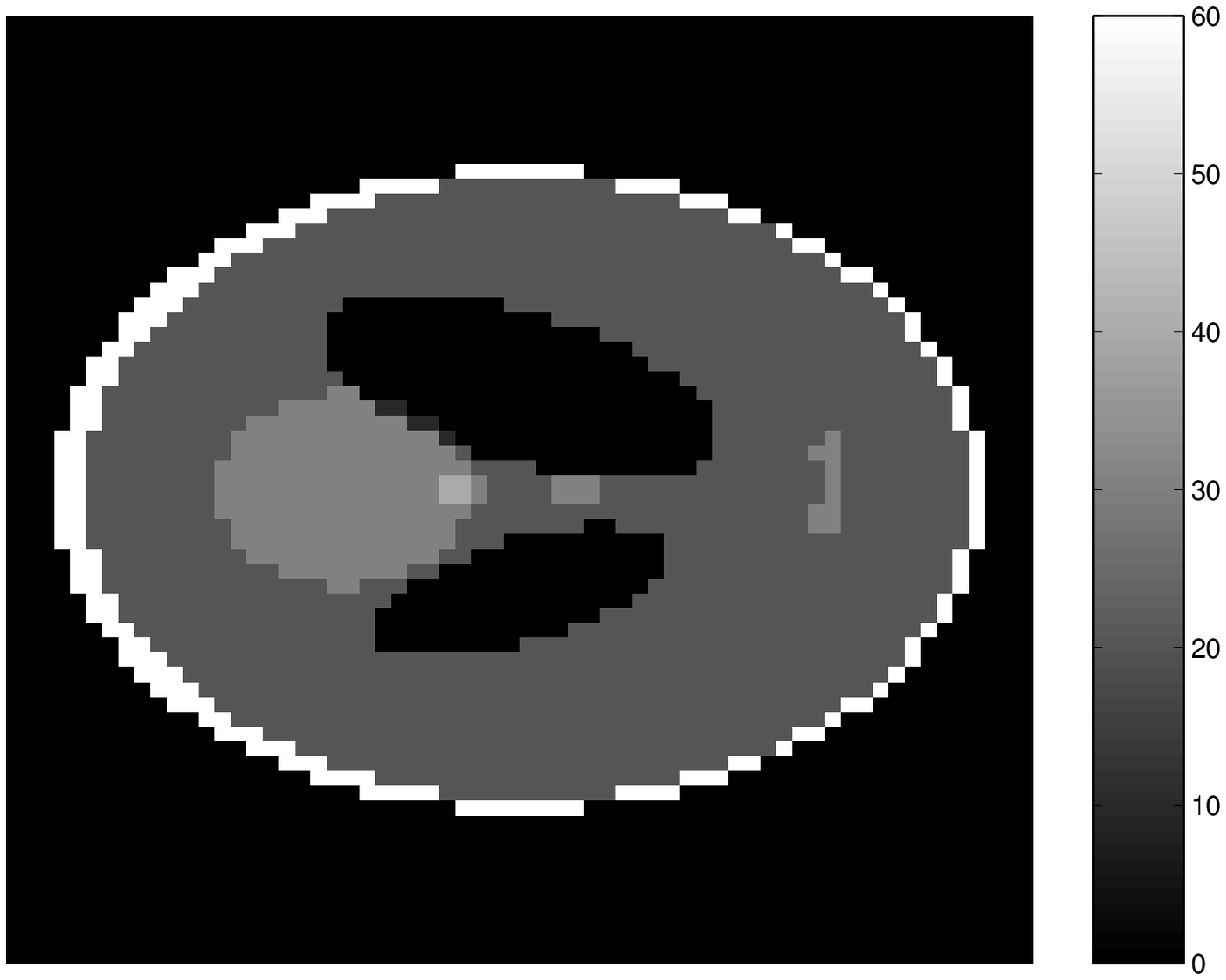}\hspace{.3cm}
\includegraphics[width=3.7cm]{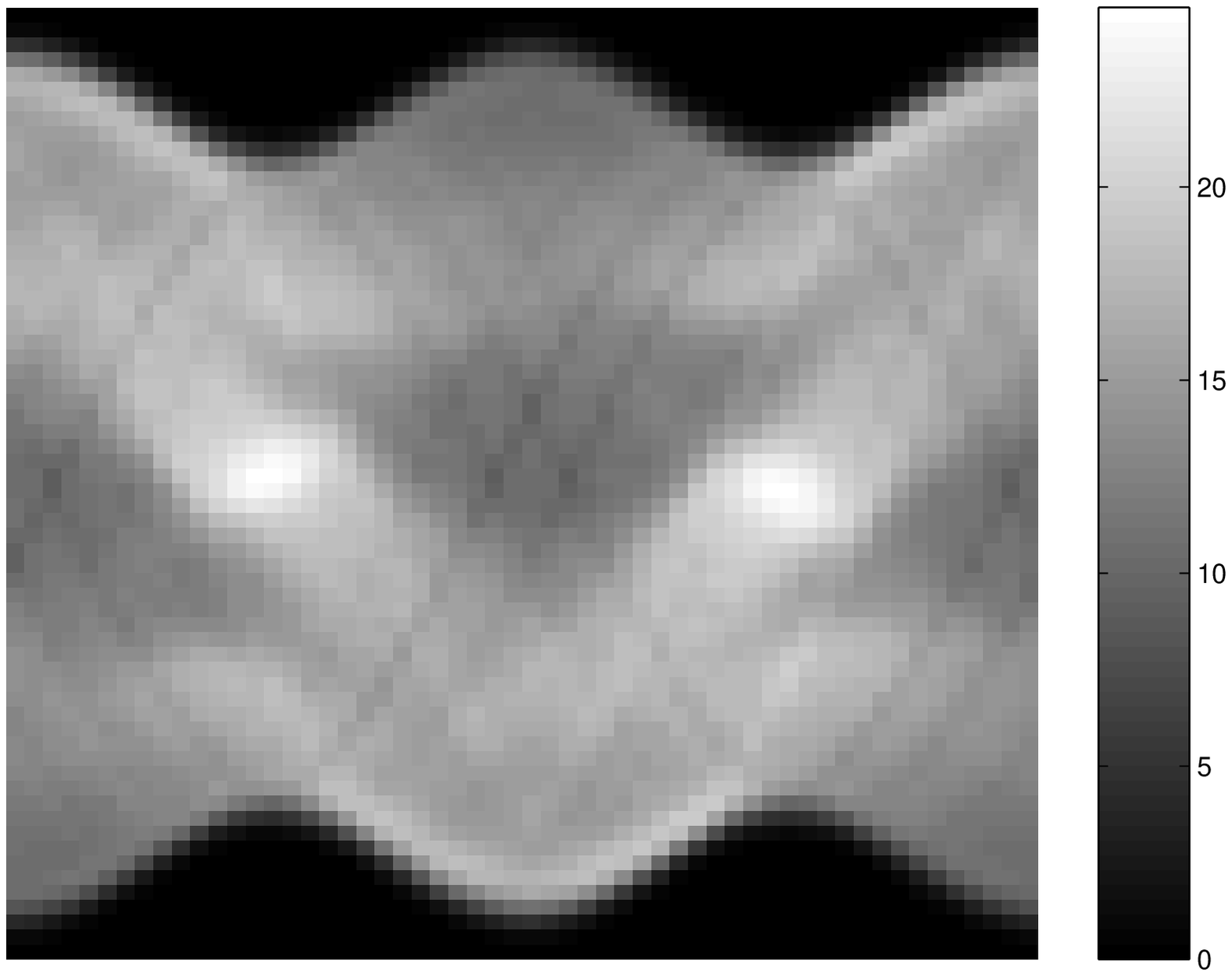}\hspace{.3cm}
\includegraphics[width=3.7cm]{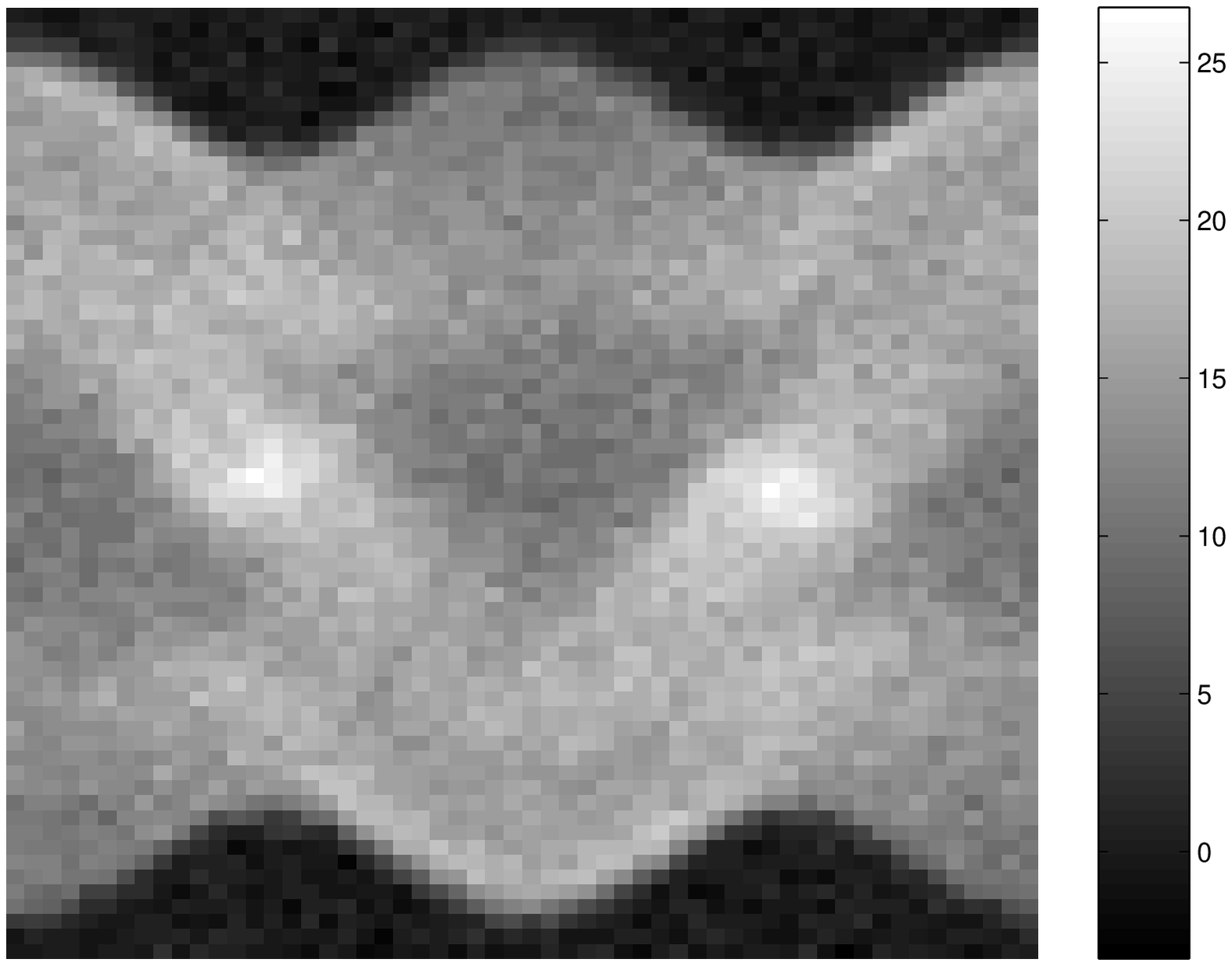}
\end{center}
\caption{\label{fig1}
The shepp-Logan phantom (left), the corresponding sinogram (middle)
and the noisy sinogram (right).
}
\end{figure}

We now illustrate our reconstruction approach by a few numerical simulation.
To begin with, we use a $64\times 64$ image of the Shepp-Logan phantom, as shown
in Figure~\ref{fig1}. The sinogram is obtained by means of an operator~$R$
which accounts for the geometrical system response,
as described in~\cite{formiconi,passeri}
(the resolution in the projection depends on the distance to the detector).
The simulated acquisition was performed with 64 bins and 64 angular values
evenly spaced over 360 degrees. An independent Poisson noise was added
to the sinogram, with parameter equal to the pixel value. (The total number
of counts in the noisy sinogram is equal to 50065.)

\begin{figure}
\begin{center}
\includegraphics[width=3.7cm]{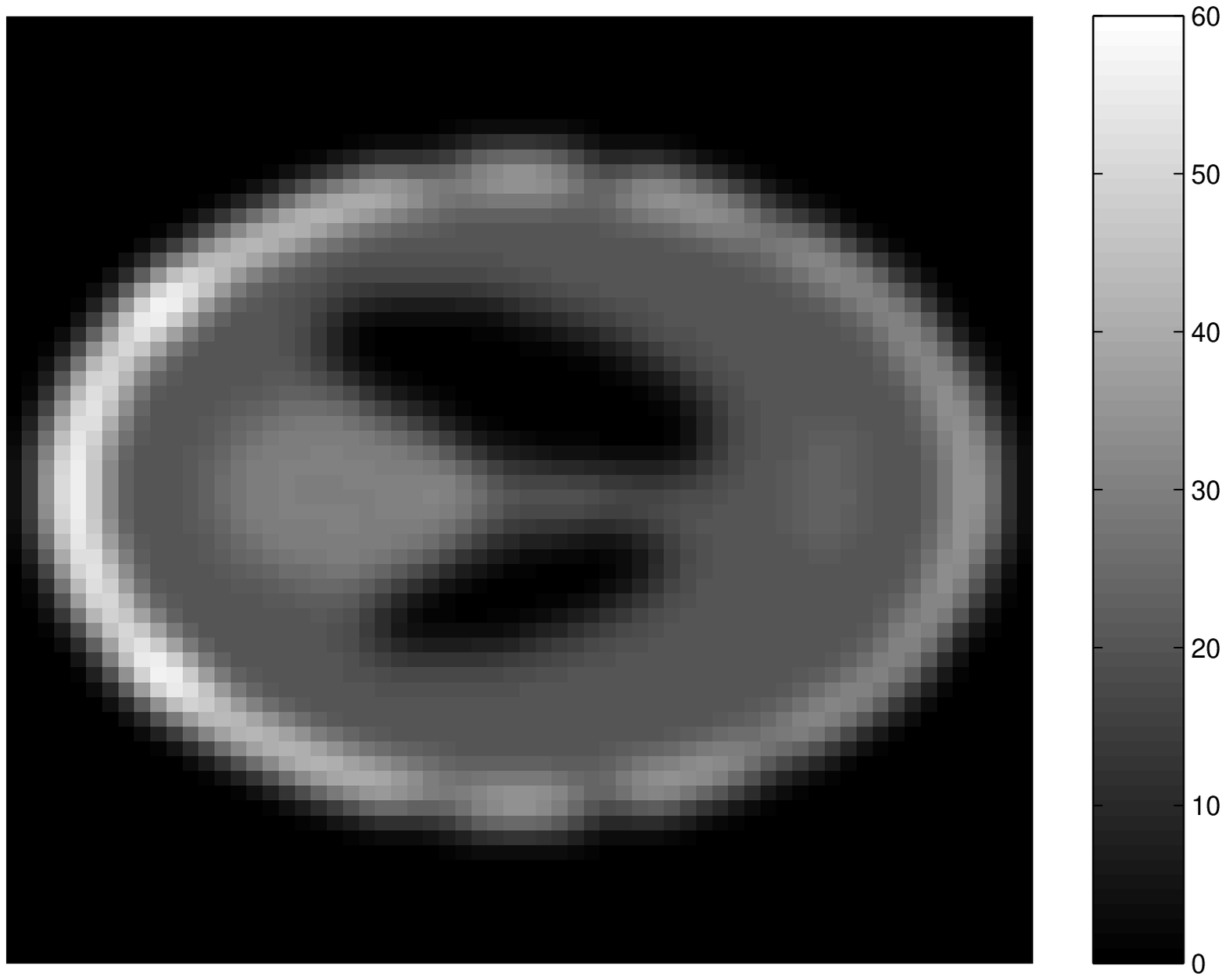}\hspace{.3cm}
\includegraphics[width=3.7cm]{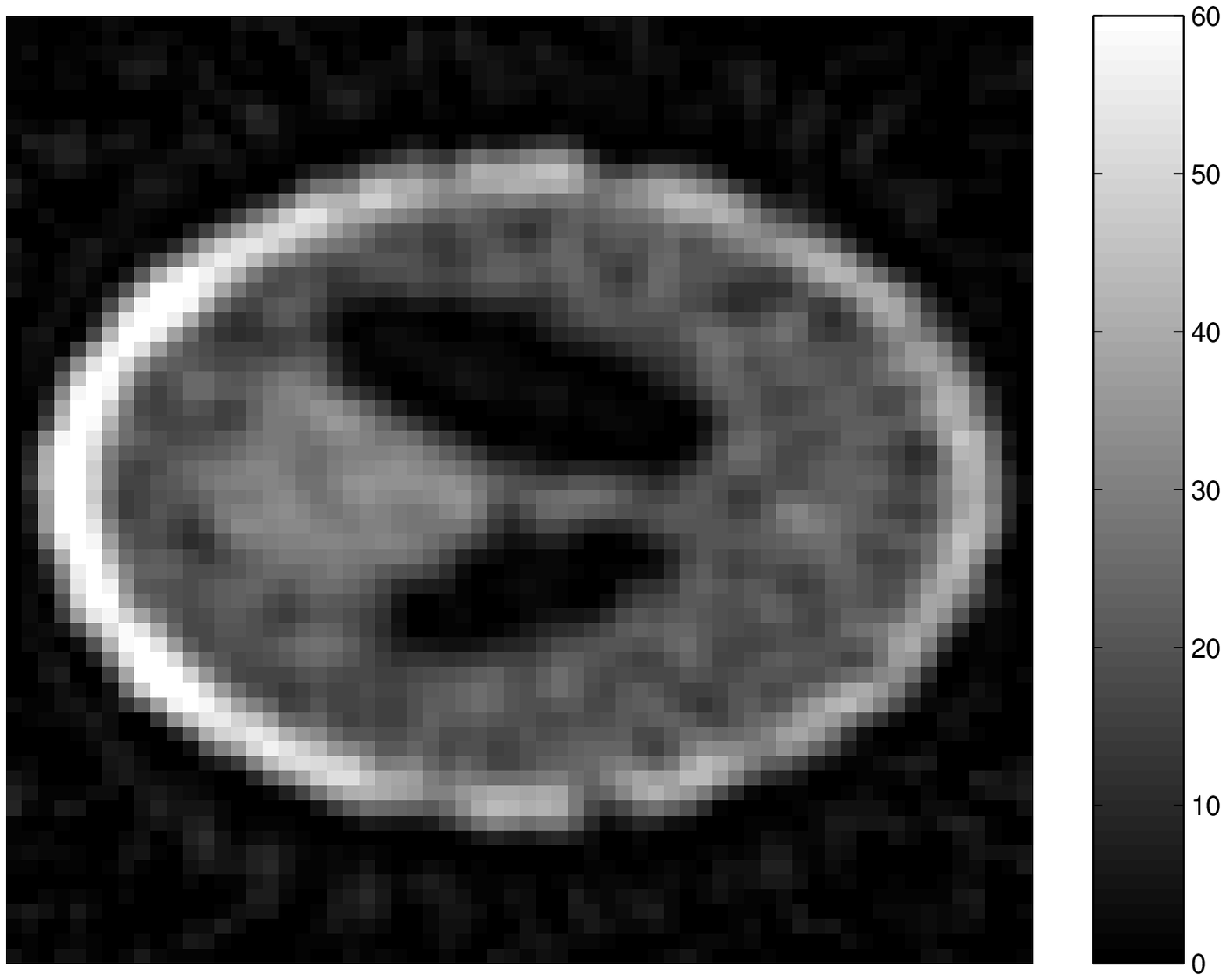}\hspace{.3cm}
\includegraphics[width=3.7cm]{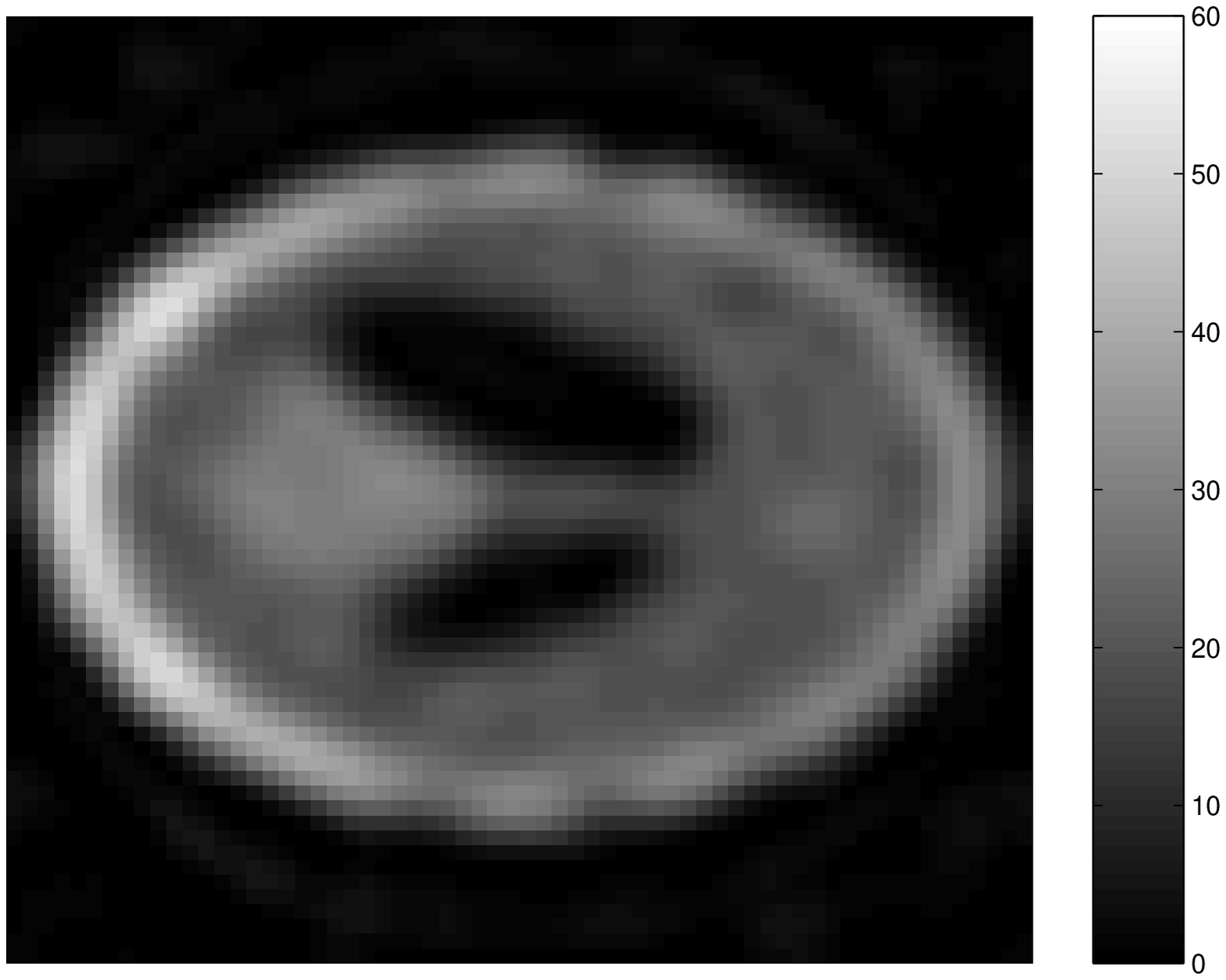}\\
\includegraphics[width=3.7cm]{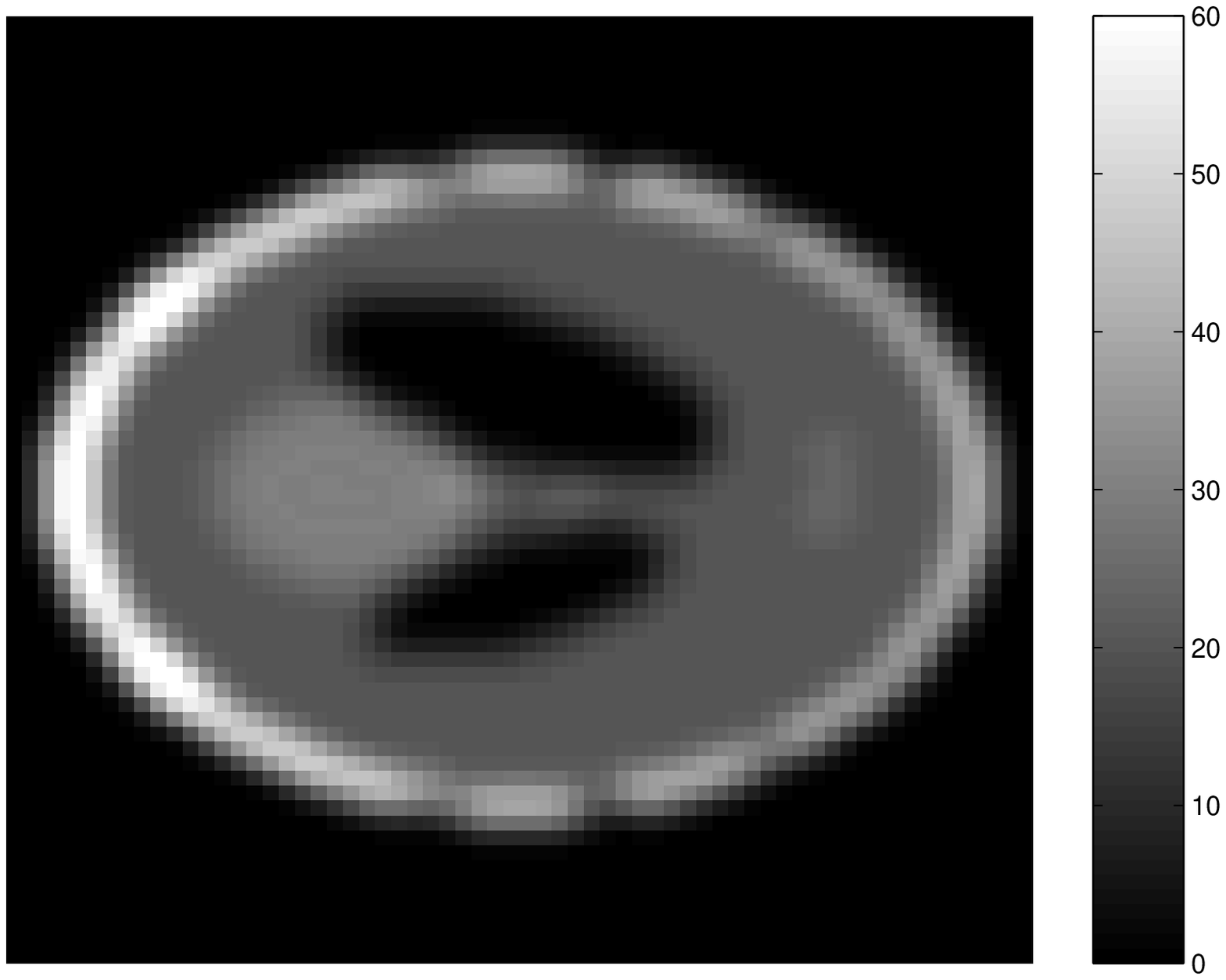}\hspace{.3cm}
\includegraphics[width=3.7cm]{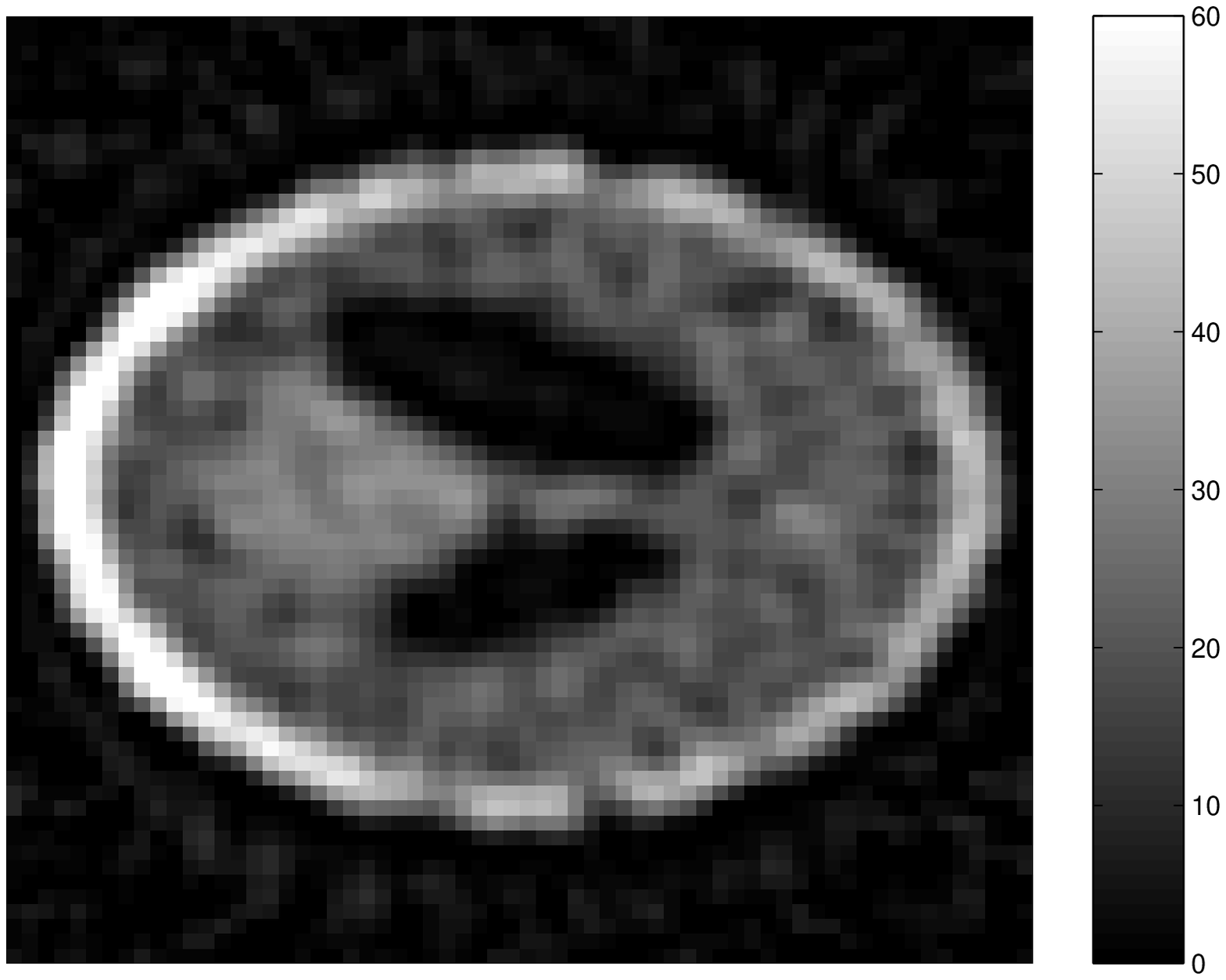}\hspace{.3cm}
\includegraphics[width=3.7cm]{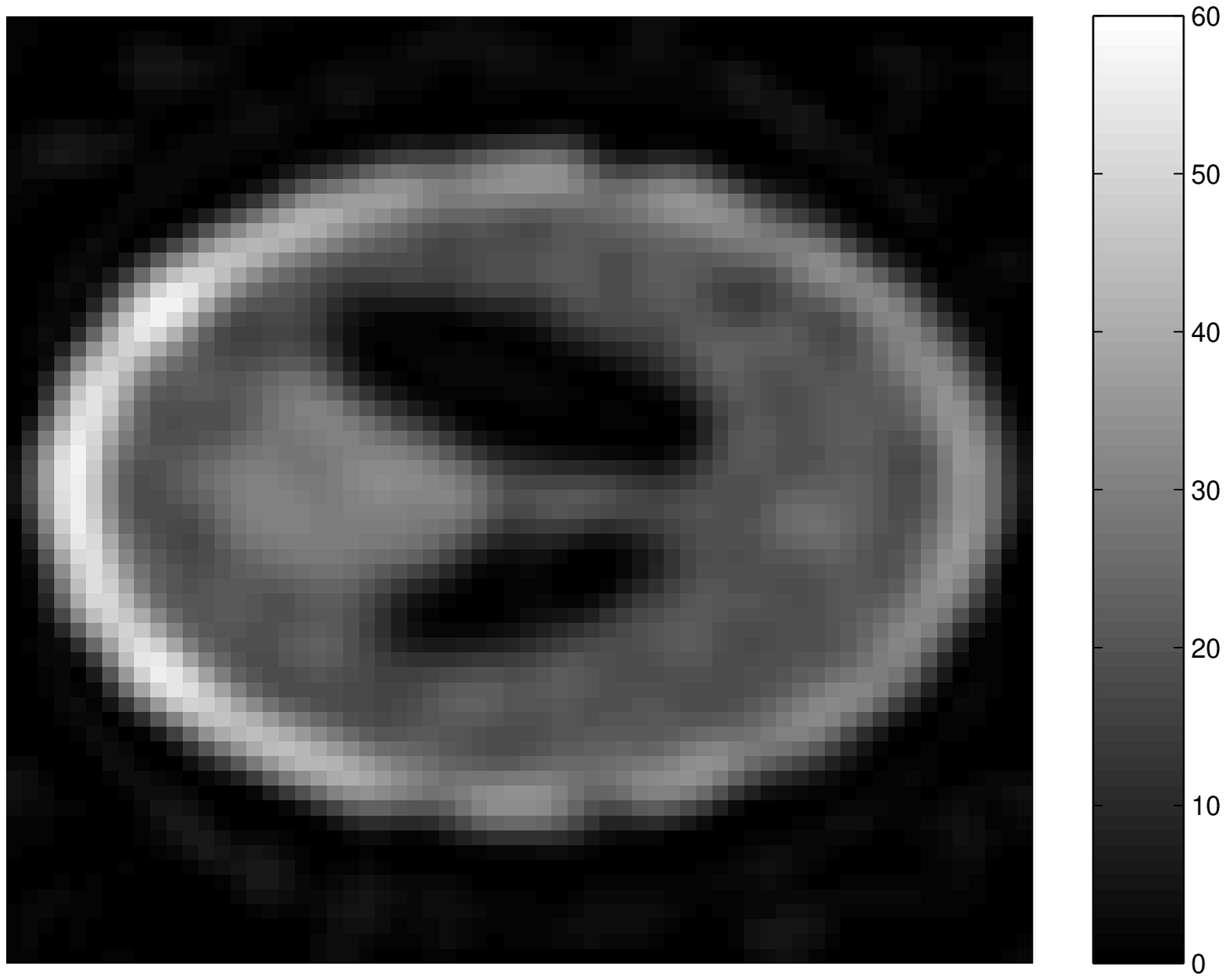}\\
\includegraphics[width=3.7cm]{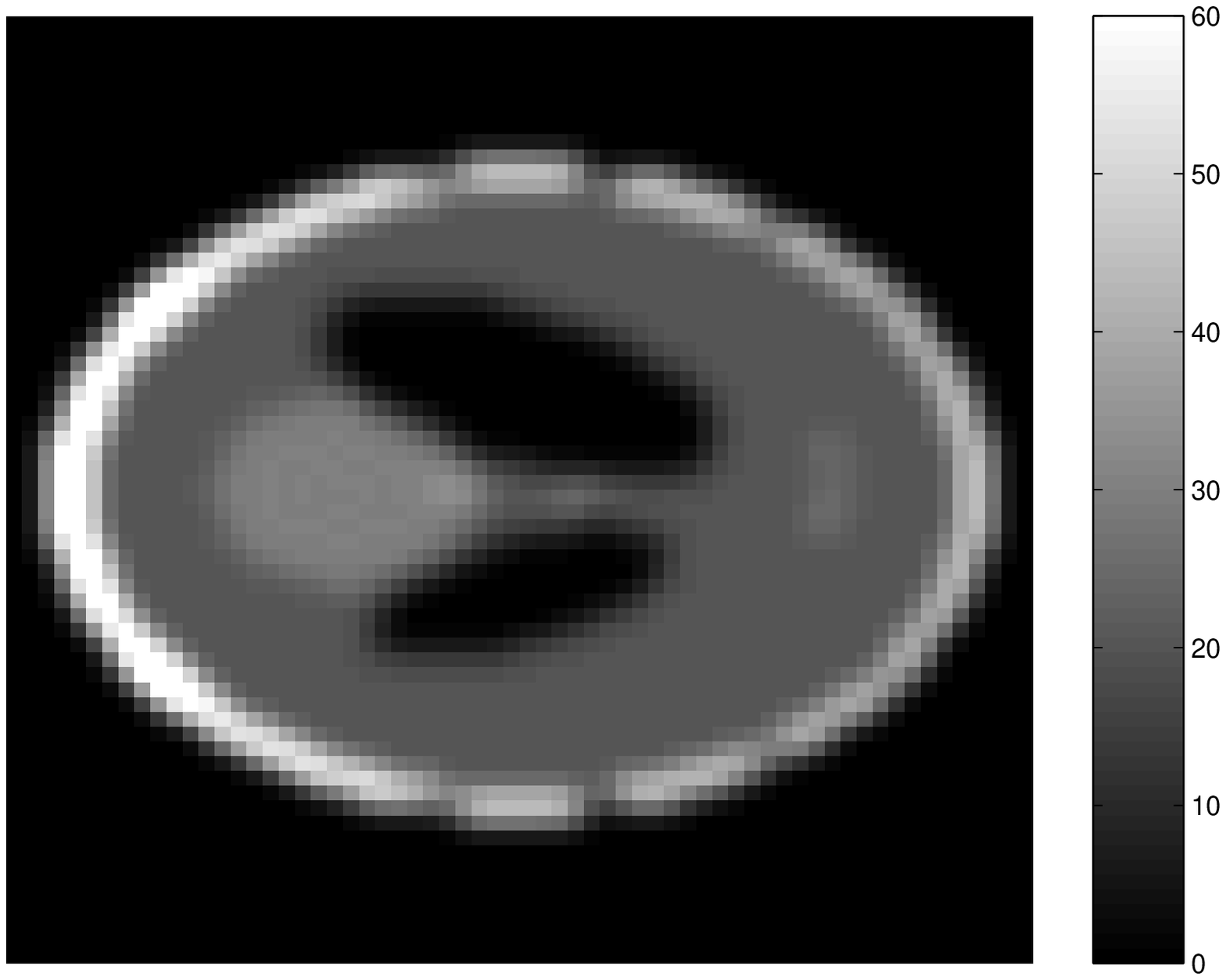}\hspace{.3cm}
\includegraphics[width=3.7cm]{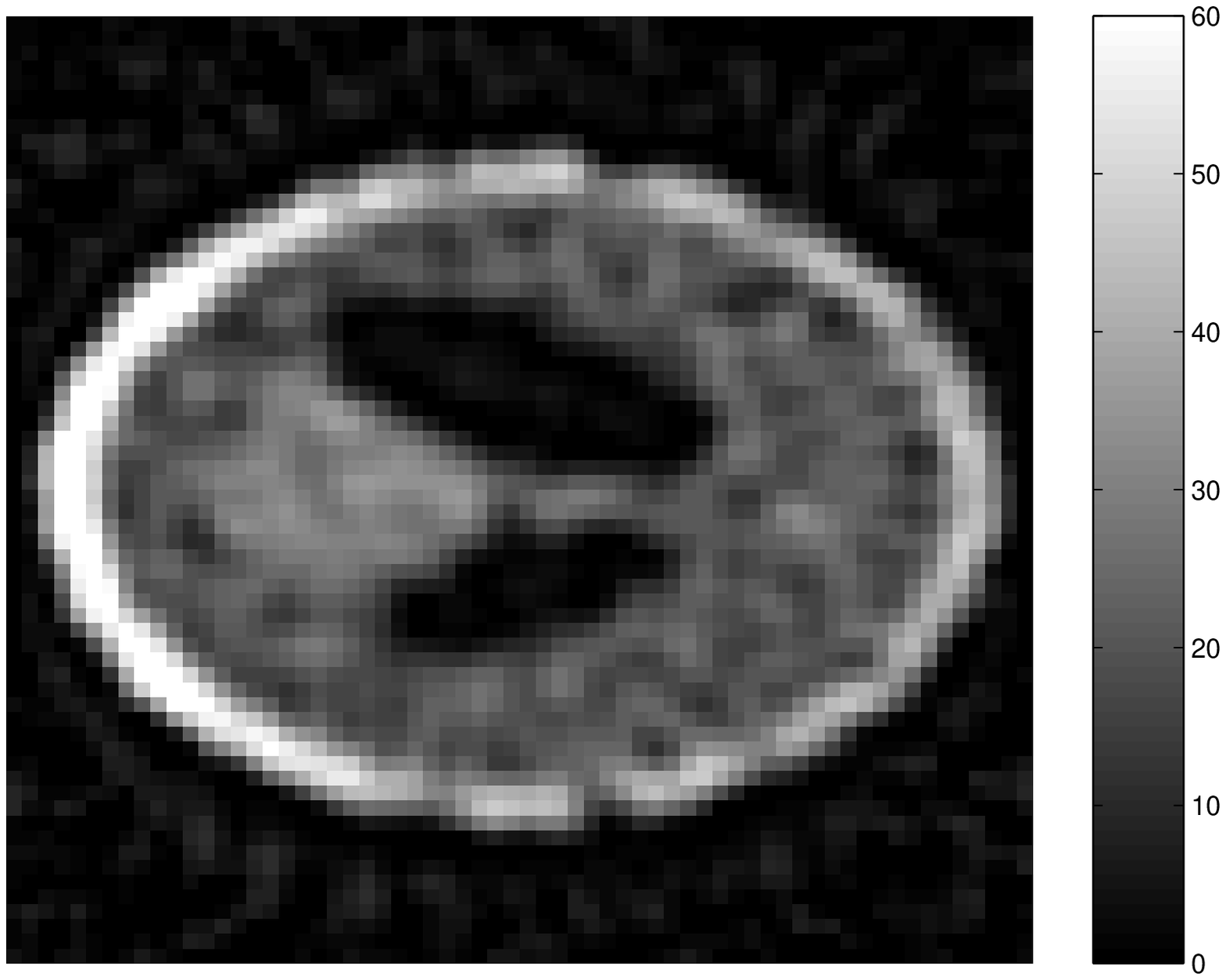}\hspace{.3cm}
\includegraphics[width=3.7cm]{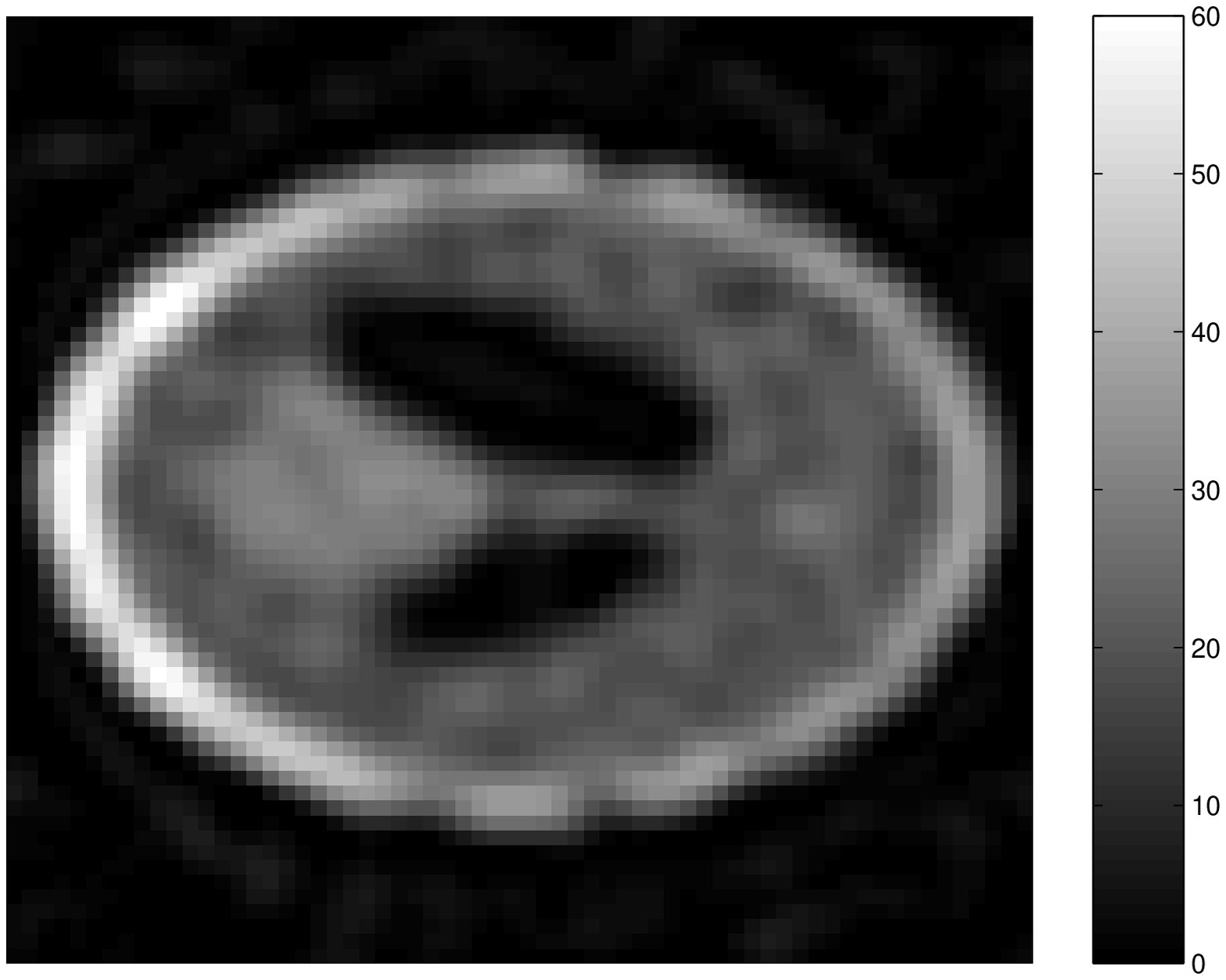}\\
\includegraphics[width=3.7cm]{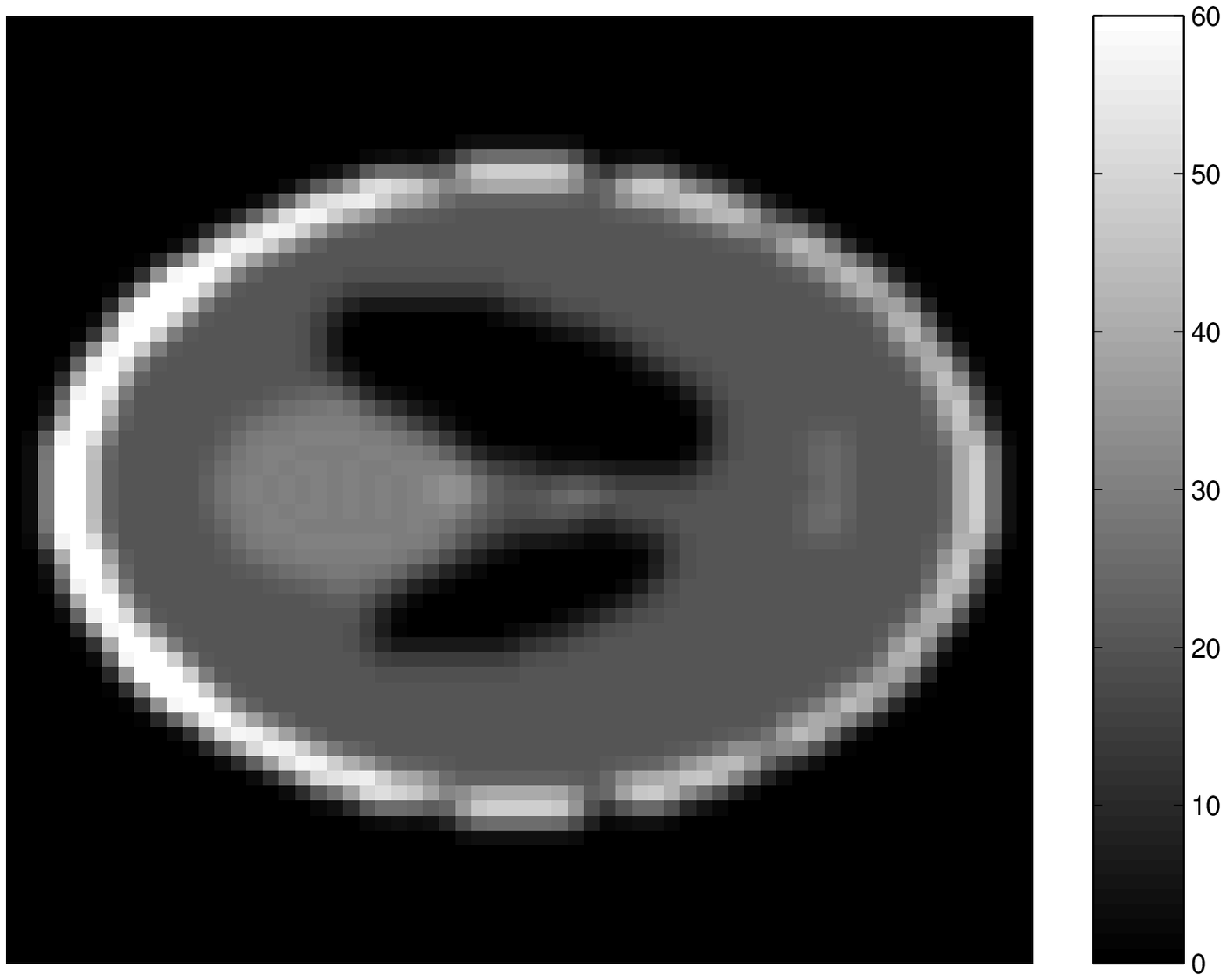}\hspace{.3cm}
\includegraphics[width=3.7cm]{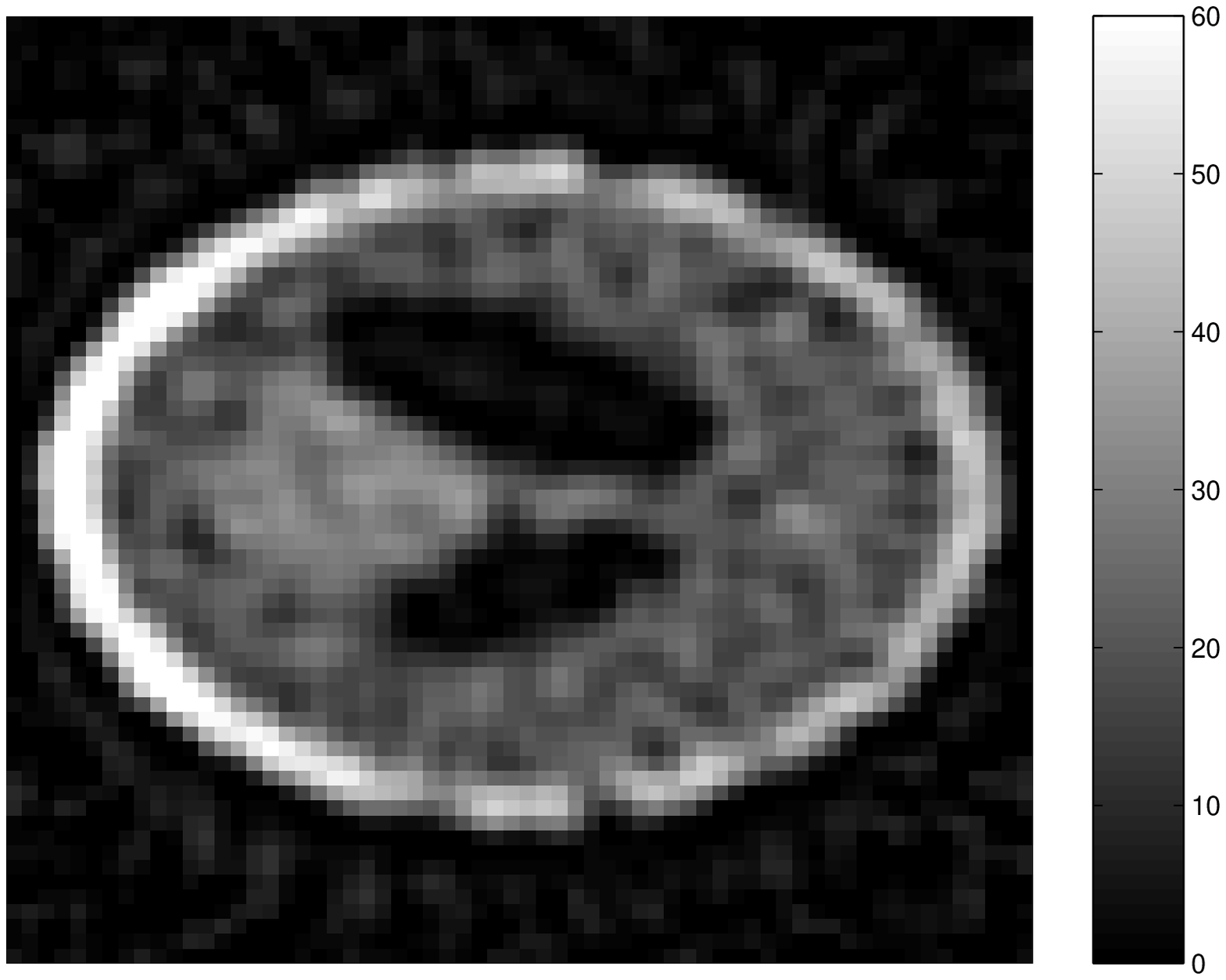}\hspace{.3cm}
\includegraphics[width=3.7cm]{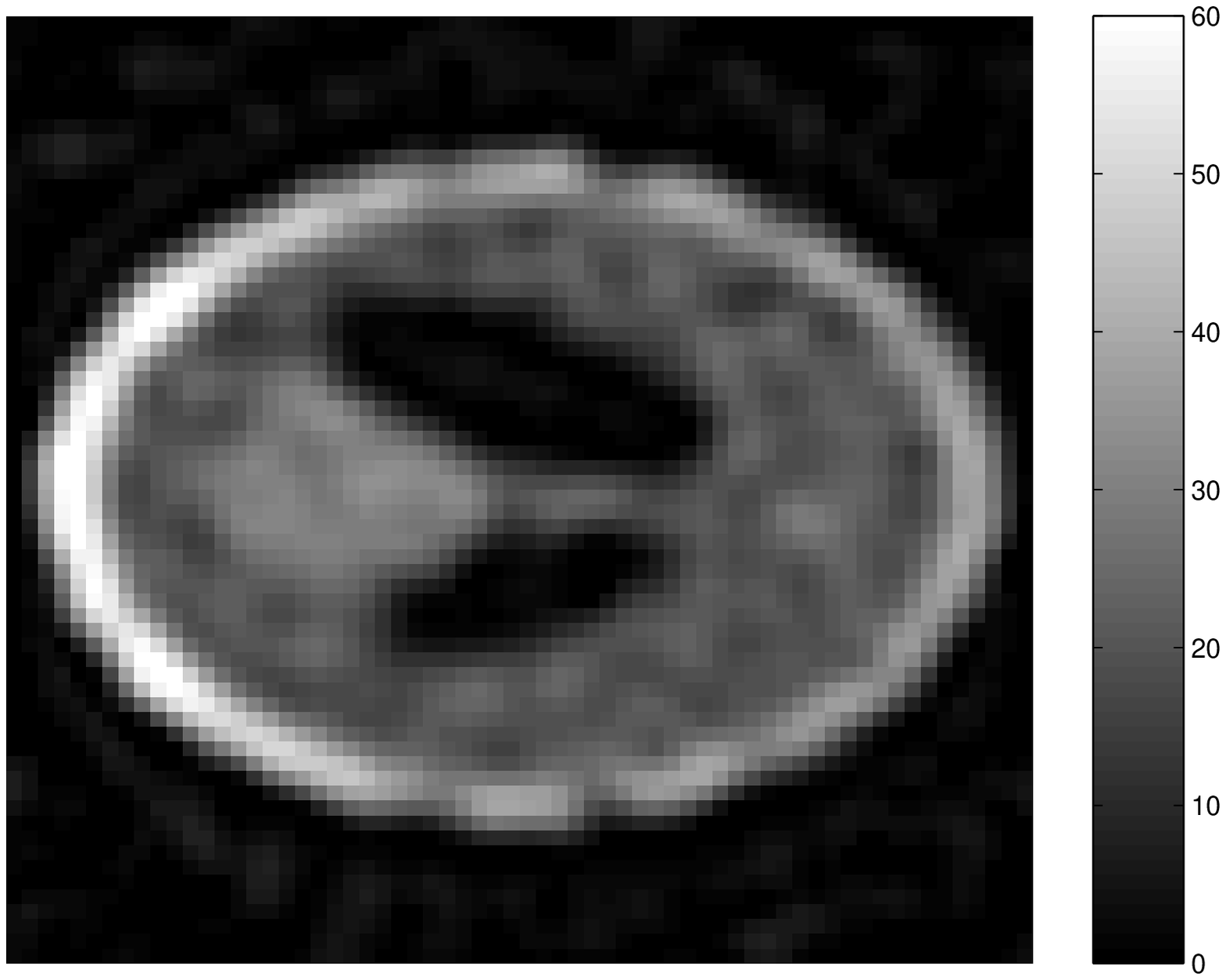}\\
\end{center}
\caption{\label{fig2}
Images in the left column are mere convolutions of the original
phantom image by the objective point-spread function (here a Hann
filter is used), with cutoff frequency equal, in Nyquist units, to
0.5 (first line), 0.6 (second line), 0.7 (third line) and 0.8 (last
line). The second and third columns show the corresponding
reconstructions, respectively without and with preprocessing
of the data.
}
\end{figure}

The reconstructions shown in Figure~\ref{fig2} were obtained by solving Problem~$(\P)$.
The data is preprocessed according to the strategy described in the
previous section. For comparison, reconstructions without preprocessing of the data
are shown.

In practice, the positivity constraint can often be
neglected: removing it from problem~$(\P)$ turns out to give
images which are essentially positive. The advantage of this is
that we deal with purely quadratic optimization, which allows for
stability analysis as well as the use of the {\sl backward error}
stopping criterion, as in~\cite{dmg-m-et-al}: the reconstructions
are performed using the conjugate gradients algorithm, with
construction of the Galerkin tridiagonal matrix.

It is clear that preprocessing enhances the quality of the reconstruction.
These simulations corroborate the anticipated relevance of our approach,
and attests its numerical feasibility. It is interesting to note that,
as conjectured in Section~\ref{computational}, the unstable step
consisting in computing $R^\dagger\vgg$ does not damage the whole
process: the application of~$C_\beta$ to $R^\dagger\vgg$ is sufficient
to ensure stability of the reconstruction. 

Table~\ref{tabNq} displays the {\sl normalized quadratic error}
$$
\E(\vgf)\eqd\frac{\norm{C_\beta\vgf_0-\vgf}}{\norm{C_\beta\vgf_0}}
$$
for various solutions~$\vgf$: in the first line, $\vgf$ is reconstructed
without preprocessing of the data, while in the second line, $\vgf$ is
obtained without preprocessing. In the last line, the value (given
for reference) is that corresponding to the reconstruction by
Filtered Back-Projection (FBP). Here, $\vgf_0$ denotes the original object.

\begin{figure}
\begin{center}
\begin{tabular}{|c||c|c|c|c|}
  \hline
&\multicolumn{4}{|c|}{Cutoff frequency} \\
\cline{2-5}
\multicolumn{1}{|c||}{ }& 0.5&0.6 & 0.7 & 0.8 \\
  \hline
without preprocessing&  0.347786& 0.301052 &  0.287429&  0.291660\\
  \hline
with preprocessing  &0.118515 &0.140563 &0.168467 &0.202067 \\
   \hline
with FBP  &0.868203 &0.879254  & 0.887703 & 0.891690\\ 
   \hline 
  \end{tabular}
\end{center} 
\caption{\label{tabNq}
Normalized quadratic error $\E(\vgf)$ for various cutoff frequencies and three
reconstruction methods: Fourier synthesis without and with preprocessing (respectively
first and second lines), and the Filtered Back-Projection (last line).}
\end{figure}

\section{Conclusion}

We have considered an extension of the regularization by mollification
of the problem of computerized tomography. This was motivated by the fact
that realistic models differ from the standard Radon transform: the latter
cannot be used if the system response accounts for geometrical aspects
as well as Compton scattering and attenuation.

The data corresponding to the objective of our reconstruction process
(namely a smoothed version of the original object) can be computed
numerically by means of the proximal point algorithm, which operates
even when the dimension of the system matrix makes it difficult or
impossible to use the SVD.

Simulations have shown both the numerical feasibility of our
regularization scheme and its efficiency, in terms of robustness
and image quality. The unstable computation of the
pseudo-inverse of the data fits in a regularization
strategy without affecting its global stability.

Such a strategy may be used to apply the regularization by mollification
to many other imaging techniques. More generally, the proximal algorithm may provide
the user with an efficient tool for preprocessing the data in accordance
with the objective of the reconstruction, whenever the latter objective
is a linear transform of the original image.

\section{Appendix: proof of Theorem~\ref{theo-conv-prox}}

We start with a finite dimensional version of Opial's lemma:
\begin{Lemma}
Let $(\vgx_k)_{k\in\Ne}$ be an $\Re^n$-valued sequence, and let $S$ be a nonempty subset
of~$\Re^n$. Suppose that
\begin{itemize}
\item[\normalfont(i)]
for all $\vgx\in S$, the sequence $(\norm{\vgx_k-\vgx})_{k\in\Ne}$ has a limit;
\item[\normalfont(ii)]
every cluster point of $(\vgx_k)$ belongs to~$S$.
\end{itemize}
Then $(\vgx_k)$ converges to a point $\overline{\vgx}$ in~$S$.
\end{Lemma}

\Proof{
Condition (i) implies that $(\vgx_k)$ is bounded, which implies in turn that
the sequence has at least one cluster point, say~$\overline{\vgx}$. We shall prove that
every cluster point~$\check{\vgx}$ must coincide with~$\overline{\vgx}$.
Let $(\vgx_{k_j})_{j\in\Ne}$ be a subsequence converging to~$\overline{\vgx}$ and
$(\vgx_{k'_j})_{j\in\Ne}$ be a subsequence converging to~$\check{\vgx}$.
Condition (ii) implies that both~$\overline{\vgx}$ and~$\check{\vgx}$ belong to~$S$.
Condition (i) again implies that
$$
\norm{\vgx_k-\overline{\vgx}}^2-\norm{\vgx_k-\check{\vgx}}^2
$$
has a limit. Developing the above expression then shows that the sequence
$\scal{\vgx_k}{\overline{\vgx}-\check{\vgx}}$
has a limit, from which we deduce that
$$
\lim_{j\to\infty}\big\langle\vgx_{k_j},\overline{\vgx}-\check{\vgx}\big\rangle=
\lim_{j\to\infty}\big\langle\vgx_{k'_j},\overline{\vgx}-\check{\vgx}\big\rangle.
$$
Therefore, $\scal{\overline{\vgx}}{\overline{\vgx}-\check{\vgx}}=\scal{\check{\vgx}}{\overline{\vgx}-\check{\vgx}}$,
which immediately yields $\overline{\vgx}=\check{\vgx}$.~\eop
}
\bigskip

Recall that the {\sl effective domain} of a convex function~$F$ is defined to be the set
$$
\dom{F}\eqd\set{\vgx\in\Re^n}{F(\vgx)<\infty},
$$
and that the {\sl subdifferential} of~$F$ at a point~$\vgx$ is the (closed convex) set
$$
\partial F(\vgx)\eqd
\set{\xgg\in\Re^n}{\forall\vgy\in\Re^n,\;F(\vgy)\geq F(\vgx)+\scal{\xgg}{\vgy-\vgx}}.
$$
Members $\xgg$ of the subdifferential are called {\sl subgradients}, and the inequality
in the definition of $\partial F(\vgx)$ is referred to as the {\sl subgradient inequality}.
\bigskip

{\sc Proof of theorem~\ref{theo-conv-prox}}.
For all~$\vgx\in\Re^n$,
$$
F(\vgx_{k+1})+\frac{1}{2\lambda_k}\norm{\vgx_{k+1}-\vgx_k}^2\leq
F(\vgx)+\frac{1}{2\lambda_k}\norm{\vgx-\vgx_k}^2.
$$
Taking $\vgx=\vgx_k$ yields
$$
F(\vgx_{k+1})+\frac{1}{2\lambda_k}\norm{\vgx_{k+1}-\vgx_k}^2\leq
F(\vgx_k).
$$
This shows that the sequence $F(\vgx_0),F(\vgx_1),\ldots$ decreases, and since~$F$
is bounded below, the sequence converges to some real limit $l\geq\eta\eqd\inf\set{F(\vgx)}{\vgx\in\Re^n}$.
In order to obtain the first assertion of the theorem, we must prove that actually $l=\eta$.
We start by proving the following inequality:
\begin{equation}
\label{inequality}
\forall\vgx\in\dom{F},\quad
\norm{\vgx_k-\vgx}^2-\norm{\vgx_{k+1}-\vgx}^2\geq
2\lambda_k\big(F(\vgx_{k+1})-F(\vgx)\big).
\end{equation}
The optimality of $\vgx_{k+1}$ reads:
$$
\vg0\in\partial F(\vgx_{k+1})+\frac{\vgx_{k+1}-\vgx_k}{\lambda_k},
$$
that is, $\vgx_k-\vgx_{k+1}\in\lambda_k\partial F(\vgx_{k+1})$. In other words, there
exists~$\xgg$ in $\partial F(\vgx_{k+1})$ such that $\vgx_k-\vgx_{k+1}=\lambda_k\xgg$.
Now, for all $\vgx\in\dom{F}$, we have:
\begin{eqnarray*}
\lefteqn{\norm{\vgx_k-\vgx}^2}\\
&=&
\norm{\vgx_k-\vgx_{k+1}}^2+\norm{\vgx_{k+1}-\vgx}^2+2\scal{\vgx_k-\vgx_{k+1}}{\vgx_{k+1}-\vgx}\\
&\geq&
\norm{\vgx_{k+1}-\vgx}^2+2\scal{\vgx_k-\vgx_{k+1}}{\vgx_{k+1}-\vgx}\\
&=&
\norm{\vgx_{k+1}-\vgx}^2+2\lambda_k\scal{\xgg}{\vgx_{k+1}-\vgx}\\
&\geq&
\norm{\vgx_{k+1}-\vgx}^2+2\lambda_k\big(
F(\vgx_{k+1})-F(\vgx)\big),
\end{eqnarray*}
in which the last inequality stems from the subgradient inequality. Thus \eqref{inequality} is clear.
Next, since $F(\vgx_{k+1})\geq l$ for all~$k$, \eqref{inequality} implies that
$$
\forall \vgx\in\dom{F},\quad
\norm{\vgx_k-\vgx}^2-\norm{\vgx_{k+1}-\vgx}^2\geq
2\lambda_k\big(l-F(\vgx)\big).
$$
Summing the above for $k=0,\ldots,n$ shows that, for all $\vgx\in\dom{F}$,
$$
\norm{\vgx_0-\vgx}^2-\norm{\vgx_{n+1}-\vgx}^2\geq
2\left(\sum_{k=0}^n\lambda_k\right)\big(l-F(\vgx)\big).
$$
Since the series $\sum\lambda_k$ is divergent, we must have $l\leq F(\vgx)$ for all $\vgx\in\dom{F}$,
which eventually show that $l\leq\eta$.

It remains to prove that, assuming $S\eqd\argmin F\not=\emptyset$, the proximal sequence $(\vgx_k)$ converges
to a point~$\overline{\vgx}\in S$.
Since $F(\vgx_{k+1})-F(\vgx)\geq 0$ for every $\vgx\in S$,
\eqref{inequality} shows that
$$
\forall \vgx\in S,\quad
\norm{\vgx_k-\vgx}^2\geq\norm{\vgx_{k+1}-\vgx}^2.
$$
Thus $\norm{\vgx_k-\vgx}$ decreases as~$k$ increases, which obviously implies
that $\norm{\vgx_k-\vgx}$ goes to a limit. It follows that the sequence $(\vgx_k)$
is bounded.
Now, by lower semi-continuity, every cluster point~$\overline{\vgx}$ of $(\vgx_k)$ (there is at least one
since the sequence is bounded) satisfies
$$
F(\overline{\vgx})\leq\lim_{j\to\infty} F(\vgx_{k_j})=\lim_{k\to\infty}F(\vgx_k)=l.
$$
Thus every cluster point of $(\vgx_k)$ belongs to~$S$, and Opial's lemma
shows that the sequence $(\vgx_k)$ actually converges to a point~$\overline{\vgx}$ in~$S$.~\eop


\begin{thebibliography}{99}

\bibitem{almasa}
		{\sc N. Alibaud}, {\sc P. Maréchal} and {\sc Y. Saesor},
		{\it A variational approach to the inversion of truncated Fourier operators},
		Inverse Problems {\bf 25}, 2009.

\bibitem{xapi}
		{\sc X. Bonnefond} and {\sc P. Mar\'echal},
		{\it A variational approach to the inversion of some compact operators},
		Pacific Journal of Optimization {\bf 5(1)}, pp. 97-110.


\bibitem{formiconi}
		{\sc A.R. Formiconi}, {\sc A. Pupi}, {\sc A. Passeri},
		{\it Compensation of spatial system response in SPECT with conjugate gradient reconstruction technique},
		Phys. Med. Biol. {\bf 34}, pp. 69-84, 1989.


\bibitem{buv98}
		{\sc P. M. Koulibaly}, {\sc I. Buvat},
		{\sc M. Pelegrini} and {\sc G. El fakhri},
		{\it Correction de la perte de résolution avec la profondeur},
		Revue de l'ACOMEN, {\bf 4}, pp. 114-120, 1998. 

\bibitem{lannes-rc}
		{\sc A. Lannes}, {\sc S. Roques} and {\sc M.-J. Casanove}, {\it
		Stabilized reconstruction in signal and image processing; Part I:
		partial deconvolution and spectral extrapolation with limited
		field}, J. Mod. Opt. {\bf 34}, pp. 161-226, 1987.

\bibitem{louis-maass}
		{\sc A.K. Louis} and {\sc P. Maass},
		{\it A mollifer method for linear operator equations of the first kind},
		Inverse Problems {\bf 6}, pp. 427-440, 1990.

\bibitem{prox-inversion}
		{\sc P. Mar\'echal}, {\sc A. Rondepierre},
		{\it A proximal approach to the inversion of ill-conditioned matrices},
		C. R. Acad. Sci. Paris, Ser. I 347, pp. 1435-1438, 2009.

\bibitem{m-togane-c}
		{\sc P. Mar\'echal}, {\sc D. Togane} and {\sc A. Celler}, {\it A
		new reconstruction methodology for computerized tomography: FRECT
		(Fourier Regularized Computed Tomography)}, IEEE, Trans. Nucl.
		Sc., {\bf 47}, pp. 1595-1601, 2000.

\bibitem{dmg-m-et-al}
		{\sc D. Mariano-Goulart}, {\sc P. Mar\'echal},
		{\sc L. Giraud}, {\sc S. Gratton} and {\sc M. Fourcade},
		{\it A priori selection of the regularization
		parameters in emission tomography by Fourier synthesis},
		Computerized Medical Imaging and Graphics, {\bf 31}, pp. 502-509, 2007. 

\bibitem{martinet}
		{\sc B. Martinet},
		{\it Régularisation d'inéquations variationelles par approximations successives},
		Revue Française d'Informatique et de Recherche Opérationelle, pp. 154-159, 1970.

\bibitem{passeri}
		{\sc A. Passeri}, {\sc A.R. Formiconi} and {\sc U. Meldolesi},
		{\it Physical modelling (geometrical system response, Compton scattering and attenuation)
		in brain SPECT using the conjugate gradient reconstruction method},
		Phys. Med. Biol., {\bf 37}, pp. 1727-1744, 1992.

\bibitem{rockafellar}
		{\sc R.T. Rockafellar},
		{\it Monotone operators and the proximal point algorithm},
		SIAM Journal on Control and Optimization, {\bf 14 (5)}, pp.877-898, 1976.

\end{thebibliography}
\end{document}